\documentstyle[nato]{crckapb}
\input psfig
\def\SBIMSMark#1#2#3{
 \font\SBF=cmss10 at 10 true pt
 \font\SBI=cmssi10 at 10 true pt
 \setbox0=\hbox{\SBF Stony Brook IMS Preprint \##1}
 \setbox2=\hbox to \wd0{\hfil \SBI #2}
 \setbox4=\hbox to \wd0{\hfil \SBI #3}
 \setbox6=\hbox to \wd0{\hss
             \vbox{\hsize=\wd0 \parskip=0pt \baselineskip=10 true pt
                   \copy0 \break%
                   \copy2 \break%
                   \copy4 \break}}
 \dimen0=\ht6   \advance\dimen0 by \vsize \advance\dimen0 by 8 true pt
                \advance\dimen0 by -\pagetotal
 \dimen2=\hsize \advance\dimen2 by .25 true in
  \openin2=publishd.tex
  \ifeof2\setbox0=\hbox to 0pt{}
  \else 
     \setbox0=\hbox to 3.1 true in{
                \vbox to \ht6{\hsize=3 true in \parskip=0pt  \noindent  
                \input publishd.tex 
                \vfill}}
  \fi
  \closein2
  \ht0=0pt \dp0=0pt
 \ht6=0pt \dp6=0pt
 \setbox8=\vbox to \dimen0{\vfill \hbox to \dimen2{\copy0 \hss \copy6}}
 \ht8=0pt \dp8=0pt \wd8=0pt
 \copy8
 \message{*** Stony Brook IMS Preprint #1, #2 ***}
}

\begin{opening}
\title{Dynamical Stability in Lagrangian Systems}
\author{Philip Boyland}
\author{Christophe Gol\'e}
\institute{University of Illinois, Urbana-Champaign and University of 
California, Santa Cruz}
\end{opening}

\def\del{\partial}
\def\ra{\rightarrow}

\def\Cross{\bigm| \kern-5.5pt \not \ \, }
\def\cross{\mid \kern-5.0pt \not \ \, }

\medskip

\def\tM{{\tilde M}}
\def\norm#1{\| #1 \|}

\def\ra{\rightarrow}
\def\de{\it}
\def\lb{\left(}
\def\rb{\right)}

\def\gm{\gamma}
\def\GM{\Gamma}
\def\lm{\lambda}

\def\om{\omega}
\def\ro{\rho}

\def\dgm{{\dot \gamma}}
\def\xd{{\dot x}}
\def\xb{{\overline x}}
\def\alx{{\alpha(x_1)}}
 
\def\implies{\Rightarrow}
\def\xo{X_c} 
\def\mo{M_c} 
\def\Mc{{\cal M}_c} 
\def\lo{S_c} 

\def\om{\omega}
\def\vmax{V_{max}}
\def\mec{{1\over 2}\norm{\xd}^2-V(x)}
\def\mect{{1\over 2}\norm{\xd}^2-V(x, t)}

\def\mml{{\cal M}_L} 
\def\mb{\overline M} 

\def\dist{dist}
\def\gm{\gamma}
\def\gmn{\gamma_N}
\def\Gm#1#2{\gm_{|_{[#1, #2]}}}

\def\GM{\Gamma}

\def\td{\tilde}
\def\ro#1#2{\delta(\gm;#1, #2)} 
 
\def\ron#1#2{\delta (\gm_N;#1, #2)} 
\def\lg{length}
\def\Lg#1#2{\lg(\gm_{|_{[{#1}, {#2}]}})}
\edef\Dist#1#2{dist(\gm ({#1}), \gm({#2}))}

\def\del{\partial}

\def\Ac#1#2{\int_{#1}^{#2}L(\gm,\dot\gm,t)dt}
\def\tae{\tilde\alpha_E}
\def\tbe{\tilde\beta_E}
\def\ae{\alpha_E}
\def\be{\beta_E}

\font\bbten=msbm10

\newfam\bbfam \textfont\bbfam=\bbten 
\def\bb{\fam\bbfam\bbten}

\def\Q{{\bb Q}}
\def\R{{\bb R}}
\def\T{{\bb T}}
\def\Z{{\bb Z}}
\def\integers{\Z}
\def\N{{\bb N}}

\def\S1{{\bb S}^1}
\def\reals{{\R}}    
\def\integers{{\Z}}
\def\P{TM\times \S1}
\def\PPP{\P}
\def\tP{T\td M\times\S1}
\def\Man{Ma\~n\'e}
\def\Mane{Ma\~n\'e}
\def\tM{\td M}
\def\Poin{{{\bb H}^n}}
\def\QED{  \rlap{$\sqcup$}$\sqcap$ }

\def\bs{{\bar\sigma}}

\def\ie{{\it i.e.  }}
\def\eg{{\it eg. }}
\def\AM{Aubry-Mather}
\def\cf{{\it cf}}

\newtheorem{theorem}{Theorem}
\newtheorem{definition}{Definition}
\newtheorem{lemma}{Lemma}
\newtheorem{proposition}{Proposition}

\newtheorem{remark}{Remark}
\newtheorem{example}{Example}

\def\Bproof{\medskip\noindent{\bf Proof\ }}
\def\Eproof{\QED\medskip}

\begin{document}

\maketitle
\SBIMSMark{1996/1b}{January 1996}{}


\section{Introduction}

The notion of stability in Dynamical Systems
refers to  dynamical behavior that persists
under perturbation.\footnote{In other contexts 
``stability'' may mean
that small perturbations of initial conditions in (some region of) 
the domain of a fixed system give rise to small changes in asymptotic behavior.
 We will always use stability to refer to persistence of the {\it global} dynamics
under perturbations of the {\it entire} system.}
By altering the nature  of the
persistence and the class of perturbations one
obtains various forms of stability.  These various forms of
stability  have proved to be extremely 
important throughout the history of dynamics.
In perhaps the best known cases, 
KAM theory and structural stability, the dynamical stability 
involves small
perturbations. However, dynamical persistence under large
perturbations (in a restricted class) is often studied and has
proved to be quite powerful. Large perturbation theories 
usually have a strong topological  component. This is  because  
behavior  that persists under  large 
perturbations must be very fundamental to the system, and 
the most fundamental aspect of a dynamical 
system is the topology of the underlying manifold.

In applying stability results, one often begins with a
model system whose dynamics are understood and then 
perturbs it.
The stability theorems indicate which dynamics of the 
model system must be present in the perturbed system.
This strategy often yields a great deal of
information about the perturbed system that could not be 
otherwise obtained.  In particular, it frequently 
provides a framework for  the investigation the other dynamics
present in the perturbed system. Since the dynamics of 
the model system must be present in all perturbed systems the
model system may also be viewed as a dynamically minimal element 
in the allowed class.

This paper surveys various results  concerning
 stability  for the dynamics of 
Lagrangian (or Hamiltonian) systems on compact manifolds.  
The main, positive results state, roughly, that if  the 
configuration manifold carries a hyperbolic metric, \ie  a 
metric of  constant negative curvature,
then the dynamics of the  geodesic flow persists in the
 Euler-Lagrange flows of a large class of time-periodic Lagrangian systems. 
This 
class contains all  time-periodic  mechanical systems on such
manifolds. 
These results are given in Theorem \ref{A} and  Theorem \ref{semiconj} 
in Section 3.  Complete proofs appear in \cite{bg}.
Many of the results on Lagrangian systems also hold for twist maps on the 
cotangent bundle of hyperbolic manifolds. 

We also present  a new stability
result for autonomous Lagrangian systems on the two torus (see Theorem \ref{t2},
Section 5) which shows, among other things, that 
there are minimizers of all rotation directions. However, in contrast
to the previously known \cite{hedlund} case of just a metric, the result allows the possibility of
 gaps in the speed spectrum of minimizers.
Our negative result  is an example of an autonomous 
{\it mechanical} Lagrangian system on the two-torus
 in which this gap actually occurs. The same system also gives us
an example of a Euler-Lagrange minimizer which is not a Jacobi minimizer on its energy level. 
 
Our results generalize several  
theories that contain what may be viewed as stability results.
The first is  the \AM\  theory. This theory  shows 
that an area preserving monotone twist map  of the annulus always 
has nicely behaved invariant sets  with each rotation number.  These 
invariant sets can be viewed 
as the remnants of the invariant circles of a minimal model, the rigid twist, 
or equivalently, the time one map of the full geodesic flow of 
the Euclidean metric on the circle. The \AM\  theory is closely 
related to Hedlund's work on geodesics on the two torus (\cite{hedlund}, \cf\ 
\cite{bangert1}).  Hedlund showed that for any Riemannian metric there 
are geodesics with all 
rotation directions, and thus the model system in this case is the 
Euclidean metric  on the torus. In related work,
closely connected to  the hyperbolic manifold results here,  
 Morse \cite{morse} showed that any metric on a 
higher genus surface  has a collection of geodesics that
``shadow'' in the universal cover the geodesics of  the hyperbolic 
metric. There are also generalizations to any dimensions and improvements
of Morse's results due to Klingenberg \cite{klingen},  Gromov \cite{gromov},
and MacKay and Denvir \cite{macden}. Because we allow
time dependent Lagrangians, these results do not imply ours.

All these theories share the property that the orbits of the 
dynamical system under consideration correspond to extremals of a 
variational problem defined in  the universal cover of the 
configuration space. The orbits that correspond to minima of the 
variational problem have special properties; they  behave 
approximately like the solutions to the  variational problem 
associated with the model 
system.  This enables one to take limits of minimizing orbit segments
or minimizing periodic orbits in order to construct a large
set of minimizing orbits on which the dynamics is similar to that of the unperturbed system. It is natural to study minimizers in the perturbed systems
because all orbits of the unperturbed system are minimizers.

There is a simple heuristic connection between these different theories.
In the Aubry-Mather
theory, Aubry's Fundamental Lemma
\cite{aubry}, \cite{meiss}, states that minimizers for a twist map
 are ordered like
orbits under a rigid rotation ( \ie  like orbits of the time-1 map of
the Euclidean geodesic flow on the circle). 
This is easily seen to imply that such orbits
have a rotation number and the rotation number of a limit of such orbits
is the limit of the rotation numbers.
If $\{x_k\}_{k\in\bf Z}$ is such
 an orbit and it has a rotation number $\om = \lim x_k/k$, then one checks that
$
\left| x_k-x_l- (k-l)\om\right|\leq 1$ for all $k, l \in \Z$. This
immediately implies that
$$
|\om||k-l| -1\leq |x_k-x_l|\leq |\om||k-l|+1.
$$
This, in Gromov's language, says that the well ordered orbits are
{\it quasi-geodesics}.  It turns out that most of the  old
or more recent Riemannian geometry results on the stability of
the hyperbolic geodesic flow can be proved using the fact that
minimal geodesics for any metric are quasi-geodesics for the hyperbolic metric. In our work on hyperbolic manifolds \cite{bg} that
we survey here, we also use a limiting argument, and our
 central Proposition \ref{qgprop}
(which is rather trivial in the autonomous case) states that {\it Euler-Lagrange minimizing segments of a
 given average speed are
quasi-geodesics.} The proof of this proposition uses techniques
that go back to Aubry's proof of his Fundamental Lemma, which
 have a  parallel in Riemannian geometry under
the guise of {\it curve shortening} arguments.
 We should remark here that the property of being a quasi-geodesic is a
 rather weak regularity property and even though it is satisfied
by all minimizers for our large class of systems on {\it any} manifolds, 
it is not sufficiently strong to make the limiting argument
work on manifolds that do not support a hyperbolic metric (except in
cases of very low dimension). 
 We illustrate this in Section 3.3 with the Hedlund
 example on the three torus.

We were greatly influenced by and used many techniques of the
recent work of Mather \cite{mather1},  \cite{mather2} (see also \cite{mane}),
 which attempts, among other things, to generalize the \AM\  theory to
higher dimensions. 
Moser showed that convex, time-periodic Lagrangian systems on compact manifolds are a natural generalization of twist maps in that twist maps
are always the time-1 maps of such Lagrangians on the
tangent space of the circle (\cite{moser}).
Denzler \cite{denzler} followed through this philosophy and gave
a proof of the \AM\  Theorem in the larger context of time-periodic, convex
Lagrangians on the circle. 
 
The class of Lagrangians we consider is almost the same as Mather's.
One fundamental difference, however, is  that Mather uses minimizers 
not in the universal cover as is done here in the hyperbolic case,
 but rather in the universal free Abelian cover
(the cover with deck group $H(M;\integers)/torsion$).  If the 
fundamental group of  
the configuration space  is torsion-free Abelian (\eg\ a torus)
this cover is the universal cover, but in general, it  is much smaller
than the universal cover.  Whereas Mather's theory works on
any compact manifold, in the special case of hyperbolic manifolds,
 much information is lost.
 Indeed, many orbits of the hyperbolic geodesic flow are not minimizers
in the sense of Mather (see Example \ref{fig8}). Hence there is no chance in finding
corresponding
orbits   in a perturbed system by looking for such minimizers.
On the other hand, in our results on the torus presented here,
we use the full strength of Mather's results (which we review
for the reader's convenience), in particular
the Lipschitz graph property of his generalized Mather sets.

It is important to remark, especially in the context of this conference,
 that Mather's stated
goal in his recent work is not so much to prove stability results 
for their own sake  but to use his Mather sets (the remnants he finds of the dynamics of the unperturbed system)  as a stairway to (generic) diffusion. He accomplished this program in \cite{mather0} for twist maps (which
are time-1 maps of time periodic Lagrangians on $\S1$), and  gets
partial results for the general case in \cite{mather2}.
 We hope that our work might help to reach this important goal.

This paper is organized as follows.
 Section 2 introduces   the Lagrangian setting and several fundamental
lemmas are given.
In Section 3, we present and outline the sketch of our two theorems 
on time periodic Lagrangians on hyperbolic manifolds.
 The first (Theorem \ref{A} in Section 3) states that you can shadow
any geodesic of the hyperbolic metric by a minimizer of the 
 Lagrangian at uniform bounded distance. This theorem also holds for symplectic twist maps in $T^*M$. The second (Theorem \ref{semiconj}) concerns large
 invariant sets for
the Euler Lagrange flow which are built  from these shadowing minimizers, and
shows that the dynamics on the sets is semiconjugate
to the hyperbolic geodesic flow.
At the end of the section we remark on  why this scheme of proof
does not work on the three torus.

In Section 4, we briefly review Mather's recent work on minimal measures,
and illustrate it by examples.
In Section 5, we apply Mather's theory to  autonomous systems on the two torus,
and give a fairly complete description of the rotation set of minimizers
 in this case (Theorem \ref{t2}, Section 5). In particular, we find
orbits of all rotation directions, and  infinitely many of
different average speeds in
each direction.  

In Section 6, we  show that the gaps in speed that are allowed
 in Theorem \ref{t2} actually occur  in natural,
mechanical systems.

\section{Preliminaries}
In this section we introduce notation and recall some basic results 
 needed in the sequel. For a thorough discussion of  Lagrangian 
systems and minimizers the reader is urged to consult Mather 
\cite{mather1}  and \Mane\  \cite{mane}. We also indicate how
to translate the different notions to the setting of
symplectic twist maps.

\subsection{Lagrangian systems.}
The main objects in the 
Lagrangian formulation of mechanics are a configuration manifold
$M$ and  a real valued function 
called a {\de Lagrangian} \index{Lagrangian function} defined on the tangent bundle $TM$.
The  configuration spaces of interest here are  closed manifolds $M$ 
with a fixed Riemannian metric $g$. The induced norm on the 
tangent bundle is denoted $\norm{v}$. We consider time-periodic 
systems determined by a  $C^2$-Lagrangian $L: TM\times \S1\to 
\R$.  The basic variational problem is to  find curves $\gm:[a,b]\to 
M$ that are extremal for the {\it action} \index{action}
$$
A(\gm)= \int_a^b L(\gm,\dot\gm,t) dt
$$
among all absolutely continuous curves $\beta:[a,b]\to  M$ 
that have the same endpoints 
$\beta(a)=\gm(a), \beta(b)= \gm(b)$.

Under appropriate hypothesis (\eg  $\gm$ is $C^1$),
such a $\gm$ satisfies
the Euler-Lagrange second order differential equations
$$
{d\over dt}{\del L\over \del v}(\gm(t),\dot \gm(t),t)-{\del L\over 
\del x}(\gm(t),\dot
\gm(t),t)=0. 
$$
Using local coordinates these equations yield a 
first order time-periodic differential equation
on $TM$, and thus in the standard way, a vector field on
$\PPP$.  Since $\PPP$ is not compact it is possible that trajectories of this 
vector field are not defined for all time in $\R$ and
thus do not fit together to give a  global flow (\ie   an $\R$-action). 
When the  flow does exist, it is called the {\de Euler-Lagrange or E-L flow}
\index{Euler-Lagrange, or E-L flow}.

We will require that Lagrangians  satisfy certain hypotheses. 
\medskip

\noindent{\bf Hypothesis:}
\medskip

{\it $L$ is a $C^2$ function $L:\PPP\ra\R$ that satisfies:
 
(a) {\it Convexity:} $\del^2 L\over \del v^2$ is positive definite.
 
(b) {\it Completeness:} The Euler-Lagrange flow determined by L 
exists.
 
(c) {\it Superquadratic:}\index{superquadratic} There exists a $C > 0$ so that $L(x,v,t)\geq 
C\norm{v}^2$.

or 

(c') {\it Superlinear}\index{superlinear}: ${L(x,v,t)\over \norm v} \to\infty$
{\it when }$\norm{v}\to +\infty$. 
}

\medskip

We will refer to $(a), (b)$ and (c) as {\it Our Hypotheses}
and $(a), (b)$ and $(c')$ as  {\it Mather's Hypothesis}. Note that
our hypotheses are a little  stronger than Mather's (adding a constant
to the Lagrangian doesn't change the E-L flow).

\begin{example} {\bf : Mechanical Lagrangians.} 
As pointed out by \Mane, Mather's Hypothesis (and hence ours) are satisfied for 
mechanical Lagrangians, \ie   those of the  form 
$$
L(x,v,t)= {1\over 2}\norm{ v}^2 -V(x,t),
$$
where the norm is taken with respect to any Riemannian metric  
on the manifold.
(In fact, one may allow the norm to vary with time, under some 
conditions, see 
\cite{mane}, page 44). 
\end{example}

\subsection{Minimizers}
Of particular interest here are extremals of the variational problem  
that minimize in the following sense.  
If $\td M$ is a cover of $M$,
 $L$ lifts to a real valued function (also called $L$)  defined on $\tP$.
A curve segment $\gm:[a,b]\to \tM$ is called a $\tM$-{\de minimizing 
segment} or an $\tM$-{\de minimizer} \index{minimizer}
if it minimizes
the action among all absolutely continuous curves $\beta:[a,b]\to 
\tM$
which have the same endpoints.

A fundamental theorem of Tonelli implies that
 if $L$ satisfies Mather's Hypotheses,
then given  $a < b$ and two distinct points
$x_a,x_b\in\tM$ there is always a minimizer $\gm$ with  $\gm(a) = 
x_a$
and $\gm(b) = x_b$. Moreover such a $\gm$ is automatically $C^2$ and
satisfies the Euler-Lagrange equations  (this uses the completeness of
the E-L flow).
 Hence its {\it differential} $d\gm(t)= (\gm(t),\dot\gm(t))$ yields
 a solution $(d\gm(t),t)$ of the E-L flow. 
A curve  $\gm:\R\ra\tM$ is called a {\de 
minimizer} if $\Gm ab$ is a minimizer for all $[a,b]\subset\R$. 
When the domain of definition of  a curve is not explicitly given it is 
assumed to be $\R$. 

As noted  in the introduction, Mather \cite{mather1} and \Man\  \cite{mane}
use $\overline M$-minimizers where
 $\overline M$ is the universal free Abelian cover. The universal cover
(which we denote $\tM$ from now on) is used here.  If $\gm$
 is an $\tM$-minimizer, we will simply say it is a minimizer.

Our main task is to get control of the speed and geometry of 
$\tM$-minimizers. 
Given a smooth curve $\gm:[c,d]\ra\tM$  and a segment $[a,b]\subset 
[c,d]$,  the {\it average
displacement} \index{average displacement} in the  cover over the interval $[a,b]$ is measured by 
$$
\ro  ab= {d (\gm(a), \gm(b))\over b-a}
$$
where $d$ is the topological metric on  $\tM$
constructed from the lift of the given Riemannian 
metric $g$. 

The fact that $L$ is superlinear or superquadratic
implies some very useful 
simple estimates on the average action of minimizers. 
These estimates are essentially in \Mane\  \cite{mane} and 
Mather \cite{mather1}. 

\begin{lemma} 
\label{disac}
Given a Lagrangian $L$ satisfying  Mather's Hypothesis,  there are functions $K\mapsto C_K^{min}, K\mapsto C_K^{max}$  both increasing to infinity such that if
$\gm$ is a minimizer and 
$ \ro ab =K$, then $$C_K^{min}K\leq {1\over b-a}{\Ac ab}  \leq 
C_K^{max}K. 
$$
In particular, if $L=\mect$, we have:
$$
{K^2\over 2}-\vmax\leq {1\over b-a}\Ac ab  \leq {K^2\over 2}-V_{min}.
$$
\end{lemma}

\Bproof
Referring the reader to \cite{mane} or \cite{bg} for a proof of the first,
general statement of the lemma, we give a proof of the mechanical case.
We first estimate the lower bound:

\begin{eqnarray*}
\lefteqn{{1\over b-a}{\Ac ab} }\\
&\geq &{1\over b-a}\int_a^b{\norm{\xd}^2\over 2}-\vmax \; dt\geq {1\over 2(b-a)^2}\lb\int_a^b\norm{\xd}\; dt\rb^2-\vmax  \\
&=&{1\over 2(b-a)^2}\lb\Lg ab\rb^2 - \vmax \geq {1\over 2(b-a)^2}\lb\Dist ab\rb^2-\vmax\\
&=&{K^2\over 2}-\vmax.
\end{eqnarray*}
The inequality on the second line follows
from  Cauchy-Schwarz:
$\int f^2\int g^2\geq \lb\int fg\rb^2$, setting $f=\norm{\xd}$ and $g=1$.

 To get the upper bound, let $\GM:[a,b]\to \td M$ with 
$\gm(a)=\GM (a)$ and $\gm (b)=\GM (b)$ be a length minimizing {\it 
geodesic} segment with respect to the given metric. Then 
$\norm{\dot \GM}= \rho (\GM;a,b)=\ro ab 
=K$. Thus, since $\gm$ is a minimizer,
$$
A (\Gm ab)\leq A(\GM)\leq\int_a^b {K^2\over 2}-V_{min}\leq \lb{K^2\over 2}-V_{min}\rb(b-a),
$$
yielding the upper bound. \Eproof


\subsection{Exact Symplectic twist maps}

For more details on symplectic twist maps, the reader is referred
to \cite{gole}
or \cite {macmeiss} (see also \cite{bernkatok} and \cite{katok}). 
An {\it exact symplectic twist map} $F$  is a map from
a subset $U$ of the cotangent bundle  of a manifold $N$ (which we
allow to be noncompact) into $U$,
which comes equipped with a {\it generating function} $S:N \times N \ra \reals$ that satisfies
$$
F^*(p\;dx)-p\;dx= P\;dX-p\;dx= dS(x, X),\eqno{(1.2)}
$$
where $(X,P)$ are the coordinates of $F(x, p)$ (this can also be
written in a coordinate free manner). 

Because the one-form
$P\;dX-p\;dx$ in (1.2) is exact, one says that $F$ is {\it exact}.  
Note that taking the exterior differential of (1.2) yields
 $dP\wedge dX= dp\wedge dx$, and so any exact $F$ is also symplectic,
\ie it preserves the standard symplectic form.
The fact that  $S$ is expressed using  the coordinates  $(x,X)$  instead of $(x,p)$ is the {\de twist condition}.  
Given $S$, one can retrieve the map (at least implicitly)
from $p=-{\del S\over \del x} and P={\del S\over \del X}$.
This can be done globally (\ie   $U=T^*N$) only
when $N$ is diffeomorphic to a fiber of $T^*N$, for example  when
$N$ is the covering space of the n-torus or of a manifold of constant negative
curvature.

The variational problem for Lagrangian systems translates into a  discrete variational problem for twist maps: the role of  curves in the continuous setting is taken by sequences of points (``integer time curves''), and the
action of a finite sequence $\xb=\{x_n, \ldots, x_m\}$
is given by $W(\xb)= \sum_n^{m-1} S(x_k,x_{k+1})$.
This corresponds closely to the continuous setting 
when the exact symplectic twist map $F$ is the time-one map of
an E-L flow. In this case,  $S(x,X)=\int_0^1L(x, \xd,t)dt$,
where $x(t)$ is the minimizer over the interval $[0,1]$
with endpoints $x$ and $X$.

In direct correspondence to  Lagrangian systems, critical points of $W$ 
(with fixed time and configuration endpoints) correspond to  orbits of $F$
(this is closely related to 
the method of  broken geodesics  in Riemannian geometry).
Action minimizers are sequences that minimize $W$ over any of their
subsegments.  The natural growth condition on the generating function 
$$
S(x, X)\geq C \; dist^2 (x, X), 
$$  
implies the  analog of Tonelli's theorem: minimizers always exist 
between any two points over any given (integer) interval of time. 
Moreover, there is 
an exact analog of Lemma 1.1: the average action of minimizers
is bounded below and above by functions of the average
displacement. The proof is virtually identical to the continuous time case, 
replacing geodesics with orbits of the time-one map of the geodesic flow.

\begin{example}
{\rm 
 Let M be $\T^n$ or a closed hyperbolic manifold, and
let $N = \tM$ be the universal cover $\R^n$ or $\Poin$,
respectively. 
On the covering space, define the  
{\de generalized standard map} using its generating function  
$\tM\times\tM\ra\reals$,
$$S(x, X)={1\over 2} dist^2(x, X) +V(x)$$
where the distance $dist$ is  induced by the Euclidean metric in $\R^n$, 
or the hyperbolic metric on $\Poin$, and $V(x)$ is $\pi_1(M)$-equivariant,
\ie  it descends to a function on $M$. A short argument shows that
one can use the relation (1.2) to solve for $(X, P)$ in terms of
$(x,p)$ and thus obtain
an exact symplectic twist
map on $T^*\tM$ that, in turn, 
 induces a map on $T^*M$ (also called a twist map).
For more general examples, \cf \ \cite{gole}.}
\end{example}

In certain cases the twist map theory overlaps with the
continuous theory. If a twist map $f$ of  $T^*\T^n$ has  a 
 generating function  that is
super quadratic in $\norm{X-x}$, the mixed  partial 
$ \partial_{12} S$ is symmetric,
and for some $a>0$ satisfies the convexity condition
$$
< \partial_{12} S (x, X).v, v >\  \leq -a \norm{v}^2
$$
uniformly in $(x, X)$, then $F$ is the  time-one map of an
E-L flow derived from a 
one-periodic Lagrangian that is superquadratic in the velocity. 
Moser \cite{moser} gives the proof in the case $n=1$. Bialy and
Polterovitch remark in \cite {bialypolt}  
that Moser's proof goes through in the case $n>1$.  This 
is not quite so, but they subsequently obtained a different
proof (personal communication).  Note that the generating function
for the generalized standard map satisfies these hypothesis.

\subsection{ Jacobi minimizers vs. E-L minimizers}

In this section, we consider autonomous mechanical Lagrangian of the form
$L(x, \dot x) = (1/2)\norm{\dot x}^2 - V(x)$ where 
$V \leq 0$ and $\norm{\dot x}$ comes from a Riemannian metric on the
manifold $M$. In this case there is a geometric
way to look for E-L minimizers on a given energy level. Recall that,
in the tangent bundle coordinates, the energy (Hamiltonian) is given by 
$H(x,\xd)= L(x,\xd)+2V(x)=
(1/2)\norm{\dot x}^2 + V(x)$ 


The norm on the tangent bundle coming from the {\de Jacobi metric} \index{Jacobi metric}
with energy $E > V_{max}$ is 
$\sqrt{E - V(x)}\norm{\dot x}$. The important fact is that geodesics
of the Jacobi metric considered as curves in $M$ 
are always the projection to $M$ of some solution curve with energy $E$
 of the Hamiltonian (or E-L) flow, and conversely, all such projections
are geodesics of the Jacobi metric. {\it However,} the parameterizations 
of the geodesics and the solutions will
usually be
different. (eg see Arnol'd \cite{arnold} or Abraham and Marsden  \cite{abrahamm})

Put in other language, this means that the extrema of the corresponding
integrals coincide in some sense. What is of importance here is whether
the  minimizers  coincide.
Let us define the notion of minimizer more carefully.

\begin{definition}

\begin{enumerate}

\item A curve $\gm:[a,b]\ra\tM$ is called an {\de E-L minimizer}
if for all absolutely continuous $\beta:[a,b]\ra\tM$ 
with $\beta(a) = \gm(a)$ and $\beta(b) = \gm(b)$,
$$
\int_a^b L(\gm, \dgm)\; dt \leq \int_a^b L(\beta, \dot\beta))\; dt
$$
\item A curve $\gm:[a,b]\ra\tM$ is called a {\de Jacobi minimizer}
for energy $E$ if for all absolutely continuous $\beta:[a_\beta,b_\beta]\ra\tM$ 
with $\beta(a_\beta) = \gm(a)$ and $\beta(b_\beta) = \gm(b)$,
$$
\int_a^b \sqrt{E - V(\gm)}\norm{\dot \gm}\; dt \leq 
\int_{a_\beta}^{b_\beta} \sqrt{E - V(\beta)}\norm{\dot \beta}\; dt $$
\end{enumerate}
\end{definition}

Note that being an E-L minimizer is a property of both the path and the
parameterization. On the other hand, being a Jacobi minimizer
just depends on the path since the integral in its definition is
$1$-homogeneous in the derivative.

The fact that solutions to the E-L equations are extrema of the 
integral in (1) is usually called Hamilton's Principle.
On the other hand, The Principal of Least Action states that if 
$\gm$ is a solution to the E-L equations 
in $\tilde M$, then among the curves $\beta$
with $\beta(a_\beta) = \gamma(a)$ and $\beta(b_\beta) = \gamma(b)$ 
which satisfy the constraint $H(\beta, \dot\beta) \equiv E$, 
the curve $\gm$ will extremize the integral $\int \norm{\dot\beta}^2$.
Now if $H(\beta, \dot\beta) \equiv E$, then $\norm{\dot\beta}^2
= \sqrt{2 (E - V(\beta))} \norm{\dot\beta}$, and so the Principle
of Least Action is the same as the variational formulation
for geodesics with the Jacobi metric as in (2).
(The reader is warned that 
there is a great deal of  variance in  the literature 
regarding the meaning of the phrases ``Principle of
Least Action'' and ``Hamilton's Principle''.)

\begin{proposition}
\label{jel} { If $\gm:[a,b]\ra\tM$ is a  Jacobi minimizer with
energy $E > V_{max}$, then there exists a re-parameterization $\gm_1$ with
$H(\gm_1, \dot\gm_1) \equiv E$  that is an E-L minimizer. }
\end{proposition}

\Bproof
Begin by noting that 
$$
L(\beta,\dot \beta) + E \geq\sqrt{2(E - V(\beta))}\norm{\dot\beta}
$$
and equality occurs if $H(\beta,\dot\beta) = E$. Indeed, setting 
$a=\sqrt{ E-V}, b=\norm{\dot \beta}/\sqrt{2}$, the above inequality is just 
$a^2+b^2\geq 2ab$, and $a=b$ means $H=E$.

 Assume that $\gamma$ is a Jacobi minimizer on the energy level
$\{ H = E\}$, and so we
may parameterize $\gm$ so it is a solution of the E-L equations.
Let this parameterization be $\gm_1:[a,b]\ra\tM$. 
Now let $\beta$ be a test path as in the definition of
E-L minimizer. Using the inequality of the previous paragraph and the
fact that $\gm$ is a Jacobi minimizer we have

\begin{eqnarray*}
\int_a^b L(\gm_1,\dgm_1) + E &=
 \int_a^b \sqrt{2(E - V(\gm_1))}\norm{\dot\gm_1}\\
 &\leq\int_a^b \sqrt{2(E - V(\beta))}\norm{\dot\beta}\\
&\leq \int_a^b L(\beta,\dot\beta) + E
\end{eqnarray*}
and so $\gm_1$ is also a E-L minimizer.\Eproof

\begin{remark} {\rm 
Somewhat surprisingly, the converse of the proposition is false;
there can be E-L minimizers which are not Jacobi minimizers on their energy
 levels.
What basically happens is that there are two E-L minimizers connecting 
the same two points
in configuration space, {\it but} they minimize over different time intervals.
They have the same energy but different Jacobi lengths, and thus they both
cannot be Jacobi minimizers. An example is given at the end of Section 6.

As noted above, The Principle of Least Action is the
same as the variational formulation of Jacobi geodesics.
Thus there are minimizers of $\int L$ with fixed time-  and
space-endpoints that are not minimizers of the action integral
with fixed space-endpoint and constant energy constraints.

Since a solution of the E-L equations is an extrema and not
necessarily a minimum of the action integral, many authors
have adapted the terminology ``Principle of Extreme Action''.                   
Our example gives another argument for this terminology: while the
extrema of the Euler Lagrange and Action problems
coincide, their minimizers may not.} 
\end{remark}

\section{Stability Results for Lagrangian Systems on Hyperbolic Manifolds} We consider here a closed hyperbolic manifold of
arbitrary
dimension $n$,
equipped with a metric of constant negative curvature. All such manifolds $M$ have the Poincar\'e $n$-disk $\Poin$ as universal covering: $\td M=\Poin$, with the canonical hyperbolic metric. Geodesics for this metric in $\Poin$ are arcs of (Euclidean) circles
perpendicular to the sphere at infinity (the Euclidean unit sphere).  The complete
proofs of Theorems 1 and 3 appear in \cite{bg}.

\subsection{Shadowing Geodesics with E-L Minimizers}
Our first result gives the  first way of formalizing the notion that Lagrangian systems satisfying our hypothesis (see Section 2.1) on hyperbolic manifolds are at 
least as complicated as the geodesics of a hyperbolic metric.
Given a hyperbolic geodesic in the Poincar\'e Disk $\Poin$, the 
theorem asserts 
that there are minimizers of the Lagrangian system that  are a 
bounded distance away and have a  variety of approximate speeds. 
Recall that $\ro ab$ means the 
average displacement in $\Poin$ over the time interval 
$[a,b]$, \ie the distance from $\gm(a)$ to $\gm(b)$ divided 
by $b-a$ (see Section 2.1).

\begin{theorem}
\label{A}
Let $(M,g)$ be a closed hyperbolic 
manifold. 
Given a Lagrangian $L$ which satisfies our hypothesis   there are 
sequences $k_i, \kappa_i, T_i$ in $\R^+$  depending only on  
$L$,  with $k_i$ increasing to infinity, such that,  for any 
hyperbolic geodesic $\Gamma \subset \Poin =\tilde M $, there are 
minimizers 
$\gm_i:\R\to\td M$ with $\dist (\gm_i, \GM) \leq \kappa_i $, 
$\gm_i(\pm \infty)= \GM(\pm \infty)$,  and 
 $k_i\leq \delta(\gm; c,d)\leq k_{i+1} $ whenever $d-c\geq T_i$.
\end{theorem}

\begin{remark}{\rm This theorem also holds for twist maps of $T^*M$ provided
their generating function has super quadratic growth, \eg the generalized standard map with generating
function $S(x, X)= (1/2) dist^2(x, X) +V(x)$. The proof is identical
to that of the Lagrangian case, using the dictionary given in
Section 2.3.}
\end{remark}

\Bproof {\bf (Sketch)}
Fix an  oriented geodesic with a given parameterization by arclength
 $\GM:\R\ra\Poin$ and a $K > K_0$ with $K_0$ as in Proposition \ref{qgprop}
below.
Let $\gmn:[-N, N]\ra\Poin$ be a minimizing segment with
$\gmn(-N) = \GM(-K N)$ and  $\gmn(N) = \GM(K N)$, and
thus $\ron{-N}{N}=K$.  We would like to be able to take
a limit of the  curves $\gmn$'s and prove that the limit is a E-L minimizer with the same endpoints at infinity as $\GM$. 
In order to do this,  we use Gromov's theory of quasi-geodesics.

\begin{definition}
\label{qgdef}
Given 
$\lambda \geq 1$ 
and
$\epsilon\geq 0$,  a curve $\gamma:\R\ra
\Poin$ or a curve segment $\gamma:[a,b]\ra\Poin$ is called
a $(\lambda,\epsilon)$-quasi-geodesic if
$$
\lambda^{-1}(d-c) - \epsilon \leq d(\gm(c),\gm(d)) \leq
\lambda (d-c) + \epsilon
$$
 for all $[c,d]$ in the domain of $\gm$.
\end{definition}

The next theorem, often called ``Stability of
quasi-geodesics'',  gives the most important property of
quasi-geodesics. It is true in the broader  context of what are
usually called $\delta$-hyperbolic spaces, but we just state the result in the 
context needed here. Given two closed subsets $X,Y\subset \Poin$,
$d(X,Y)$ denotes their Hausdorff distance  as induced by the 
hyperbolic
metric. For a proof and more information see \cite{ghys}:
\begin{theorem}
\label{qgthm}
Given  $\lambda \geq 1$ and $\epsilon \geq 0$, 
there
exists a $\kappa>0$ so that whenever $\gm$ is a $(\lambda, 
\epsilon)$-quasi-geodesic segment  in $\Poin$ and $\Gamma_0$ is 
the geodesic 
segment connecting
the endpoints of $\gm$, then $d(\gamma, \Gamma_0) <\kappa$.
If $\gm$ is a $(\lambda, \epsilon)$-quasi-geodesic, then 
$\gm(\infty)$ and $\gm(-\infty)$ exist and further, if
$\Gamma$ is the geodesic connecting
$\gm(\infty)$ and $\gm(-\infty)$, then 
 $d(\gamma, \Gamma) < \kappa$.
\end{theorem}

\begin{proposition} 
\label{qgprop}
There exists $K_0 > 0$ depending only on $L$ and
 $\lambda > 1$ and $\epsilon>0$ 
depending only on $K$ and  $L$  so that 
whenever $\gm:[a,b] \to \td M$ is a minimizing segment with $\ro ab  
=K > K_0$ and $b-a\in\N$, then 
 $\gm$ is a $(\lambda, \epsilon)$-quasi-geodesic segment.
\end{proposition}

Thus, our proposition implies that the $\gmn$ defined in the first paragraph
of the proof of Theorem 1
 are $(\lm,\epsilon)$-quasi-geodesics for the {\it same} $(\lm,\epsilon)$. Theorem \ref{qgthm} implies that the $\gmn$'s stay at uniformly bounded
distance from $\GM$.
This enables us to  take  a (pointwise) limit of the $\gmn$'s and 
get an E-L minimizer $\gm$ which has same endpoints at infinity as $\GM$ and
stays at bounded distance from $\GM$. 

The proof of  Proposition \ref{qgprop} is somewhat technical.
We just sketch the essential features.
We have that our minimizer $\gm$ satisfies
$\ro ab =K$. Since $\ro ab= {\Dist ab\over b-a}$, this  translates to $ \Dist ab = K(b-a)$. From this estimate at the endpoints of $\gm$,
we want to derive similar estimates for {\it any} interval $[c,d]\subset [a,b]$. Ideally, we would like to prove that
\begin{eqnarray}
\label{K'K''}
K'\leq \ro cd \leq K''
\end{eqnarray}
for $K', K''$ only depending on $K$. This 
would immediately prove  Proposition \ref{qgprop}. The second inequality
follows from a results of Mather:
$\ro ab =K\implies \norm{\dot \gm}\leq K''$ where $K''$ only
depends on $K$. The first inequality, proved in \cite{bg},
is only true when $(d-c)>N$, for some $N$ only depending on
$K$. But this restriction only adds an  $\epsilon$ in 
the quasi-geodesic estimate. Note that, {\it in the autonomous case}, these
estimates are trivial.

To prove the second inequality in Formula \ref{K'K''},
 Mather uses a surgery argument of 
a type 
which is widespread in the Aubry-Mather theory. If $\gm$ was going
very fast on a subinterval $[c,d]$, then one could cut $\Gm cd$
out and replace it by a minimizer between the same
endpoints but on a greater interval of time, thus reducing
the average speed on this interval. 
To make up for lost time one does the same kind of surgery
on another judiciously chosen interval on the curve. The action of this
new curve $\gm^*:[a,b]\to \td M$ is then estimated using 
the average action vs. average speed estimate of Lemma \ref{disac}
of section 2.2, and shown to be  less
than that of $\gm$, which is absurd. Our proof of the first inequality
uses similar techniques.

To finish the sketch of the proof of Theorem \ref{A}, we note that,
in our construction of the minimizers $\gmn$'s, we can vary $K$. Because of the 
bounds on the average speed on subintervals we explained above,
 this implies that, if $K_i$ is chosen to be 
increasing sufficiently fast to $\infty$ with $i$, the corresponding limiting $\gm_i$ 
are distinct and of average speed increasing to infinity.

\subsection{ Semiconjugacy with the Geodesic Flow}
One way to formulate the fact that the perturbed system 
is at least as complicated as the model system (\ie the dynamics of
the model systems don't go away) is to show that the perturbed 
system always has a invariant set that carries the dynamics of  the 
minimal model. More precisely, one shows that there is a
compact invariant set that is semiconjugate to the minimal model 
(this strategy is common in more topological theories, see
\cite{philtopo}).
MacKay  and Denvir \cite{macden} have recently extended Morse's results to the 
case with boundary and proved a result giving this semiconjugacy. 
Also Gromov \cite{gromov} and others have done this in the case of geodesic 
flows.

\begin{theorem}
\label{semiconj}
Let $M$ be a closed hyperbolic manifold 
with a 
hyperbolic metric $g$ with geodesic flow $g_t$. Given a
Lagrangian $L$ which satisfies the Hypotheses of Section 2  with
E-L flow $\phi_t$, there exists a sequences $k_i$  and $T_i$  
with $k_i$ increasing to infinity, and a family  of compact, 
$\phi_t$-invariant  sets $X_i \subset \PPP$  so that
for all $i$, $(X_i, \phi_t)$ is semiconjugate to $(T_1 M, g_t)$ 
and  $k_i\leq \delta(\phi_t(x);0,T)\leq k_{i+1}, $
whenever $T\geq T_i$ and $x\in X_i$.
\end{theorem}

Note that the geodesic flow of a hyperbolic metric is 
transitive Anosov and thus  is Bernoulli, 
has positive entropy, etc. Thus Theorem B implies that the E-L flow is 
always dynamically very complicated.

We first state precisely
what we mean by {\de semiconjugacy} \index{semiconjugacy}: 
\begin{definition}
Two flows $(X,\phi_t)$ and $(Y,\psi_t)$ are said to be 
{\it semiconjugate} (or sometimes orbit semi-equivalent) if there
is a continuous surjection $f:X\ra Y$ that takes orbits of
$\phi_t$ to those  of $\psi_t$ preserving the direction of the flow, 
but not 
necessarily the time parameterization. Note that $f$ is a local
homeomorphism  when 
restricted to an orbit of $\phi_t$, but $f$ may take many orbits of 
$\phi_t$ to the same orbit of $\psi_t$.
\end{definition}

Given a $K$ and a geodesic $\GM$,  we have constructed a $(\lm,\epsilon)$-quasi-geodesic $\gm _\GM^K$ which
shadows $\GM$.  The differential of such a curve is in $T\tM$. Project all these differential curves down to $TM$ and take the closure of this set  in $TM$.  This  gives a compact invariant set $Q_K$ for the E-L flow because the velocities of the $\gm_\GM^K$ are uniformly bounded. Finally, define $\td Q_K$ to be the set of all possible lifts to $T\tM$ of all the points in $Q_K$ .

We will show that the E-L flow restricted to $Q_K$ is semiconjugate to the
geodesic flow. $\td Q_K$ is a
 set made of $(\lm, \epsilon)$-quasi-geodesics,
for some fixed $(\lm,\epsilon)$. Each such quasi-geodesic 
shadows a unique geodesic, by Theorem \ref{qgthm}, and each
geodesic is shadowed by  at least one
 quasi-geodesic in $\td Q_K$, by
Theorem \ref{A}. Hence we have a well defined application $z\mapsto \GM_z$ which, to a point $z$ in $\td Q_K$, makes correspond the unique geodesic $\GM_z\subset \Poin$ that the
quasi-geodesic  to which $z$ belongs shadows.

We then project $z$ on $d\GM_z$ in the following fashion.
Fix a parameterization by arclength for each geodesic $\GM$.
Project $\pi(z)\in \Poin$ on $\GM_z$ via the orthogonal
projection (\ie  by drawing the unique geodesic through
$\pi(z)$ which is perpendicular to $\GM_z$) and get a point
$\GM_z(s(z))$ on $\GM_z$.  Define
$$
\sigma(z)= \lb\GM_z(s(z)), \dot \GM_z(s(z)\rb.
$$
 Theorem \ref{A}  implies that $\sigma: \td Q_K\to T_1\tM$ is onto. It is not too hard to see that $\sigma$ is also continuous,
and equivariant,  \ie it descends to a continuous map $Q_K\to T_1M$. Unfortunately, it is not necessarily injective when restricted
to an orbit of the E-L flow in $Q_K$. This is remedied using an averaging
technique due to Fuller \cite{fuller}.

Given $z=z(0)\in Q_K$ and its orbit $\{z(t), t\in \R\}$ under the E-L flow, let $a(z,t) = 
s(z(t))
- s(z(0))$, 
where,  $s(z)$ is as above.  If $a(z,t)$ were positive
for all $t$, we would be done: $\sigma$ would be injective along 
$z(t)$. Because we cannot assume that, we use the fact that $a(z,t)$
is positive for large $t$ and ``average'' $\sigma$ over the interval
$[0,t]$. 
Given ${\alpha} > 0$ define 
$$
\bar\sigma_{\alpha}(z) = \sigma(z) + 
{1\over{\alpha}}\int_0^{\alpha} a(z,t)\; dt.
$$
Informally,  $\bar\sigma_{\alpha}(z)$ is the average value of 
$\sigma$ over
 the orbit segment $z([0,{\alpha}])$. 

Now 
since for every $z\in Q_K$ we have that $\omega(z) = 
\GM_z(\infty)$, it
follows that for each $z$ there is an $\alpha_z$ so that $a(z,\alpha_z)
> 0$. Since $Q_K$ is compact, we may find an $\alpha$ with
 $a(z, \alpha) > 0$ for all $z\in Q_K$. Let $\bs =
\bar\sigma_{\alpha}$. As before, we write $\bs(z)= (\GM_z(\overline s(z)) , \dot \GM_z(\overline s(z))$. Now $\bs$ is clearly continuous, equivariant,
onto and takes orbits to orbits. We then show \cite {bg} that it is injective on
orbits $z(t)$ in $Q_K$ by showing that for any $t>0$, $\overline s(z(t))-\overline s(z(0))> 0$. This is done by using the
fact that $a$ is an additive cocycle,  \ie $a(z, t_1 + t_2) =
a(z, t_1) + a(\tilde\phi_{t_1}(z), t_2)$, for all $t_1,t_2$. By repeating
this whole process by setting $K=K_i$, where $K_i$ is as in Theorem
\ref{A}, we get distinct $Q_{K_i}$ which are semiconjugate to the
geodesic flow, and such that the average speed of a E-L orbit on
$Q_{K_i}$ goes to $\infty$ with $i$. This finishes the sketch of our
proof of Theorem \ref{semiconj}.

\subsection{Why the Limit Argument Does Not Work in $\T^3$}

Quasi-geodesics can obviously be defined on any manifold (in
fact, on any metric space).
Proposition \ref{qgprop} is valid on any compact manifold,
not just those that  support a hyperbolic metric. The reason
 why the proof above does not apply to $M=\T^3$  must therefore
be that  quasi-geodesics on $\R^2$ do not satisfy Theorem \ref{qgthm}.
Thus a sequence of  $(\lm,\epsilon)$-quasi-geodesic segments with
endpoints on a fixed geodesic $\GM$ do not have 
to stay at bounded distance of $\GM$. 

As a very simple example,
one can look at $\R^2$ with the
Euclidean metric and take the geodesic $\GM$ to be the line $y=x$. Take
the sequence $\gmn$ of ``corner'' curves going between  the
points $(-N,-N)$ and
$(N,N)$ of $\GM$ by first following the line $x=-N$ upward at unit speed
until it reaches $y=N$, and then follow that line  to $(N,N)$.
It is easy to check that the $\gmn$'s are quasi-geodesics segment
for a 
$(\lm,\epsilon)$ uniform in $N$. But obviously $d(\gmn,\GM)\to 
\infty$. 
It is not too hard either to see that one can adapt this argument in
Hedlund's example (see Section 4.2) in $\T^3$ to show how our 
limiting scheme would fail in finding a minimizer  of
rotation vector $(1,1,0)$ there: in this case, the minimizers $\gmn$'s  between
the points $(-N,-N,0)$ and $(N,N,0)$ of the (Euclidean) geodesic $\GM$ 
given by the line $y=x, z=0$
of $\R^3$  would
be corner curves that follow the tubes, jumping tubes at the
beginning, corner and end of the curve. Again the $\gmn$ are quasi-geodesic, with uniform $(\lm,\epsilon)$
by Proposition \ref{qgprop}, and again their distance to $\GM$ goes
to infinity. Note that $\gmn\to\GM_x\cup \GM_y$ in the Hausdorff
topology,  where $\GM_{x }, \GM_ y$ are the periodic orbits projection of the tubes parallel to the $x$ and $y$ axes. This is the Mather set for the vector of
rotation $(1,1,0)$ (see next section).

A very similar phenomenon happens in Example \ref{fig8} of Section 4.2. In the case
considered there the $\gamma_N$ with ``corners'' show up in the universal 
free Abelian cover $\mb$.

\section {A Quick Review of Mather's Theory of Minimal Measures}

\subsection{Theory}

 For a more detailed exposition the reader is urged to consult Mather \cite{mather1} or Ma\~ne \cite{mane}. Mather also gives a very nice survey of this theory in the beginning of \cite{mather2}.
 Given a E-L invariant probability measure 
with compact support $\mu$ on $TM\times \S1$, one can define its 
{\it rotation vector}\index{rotation vector}
 $\rho (\mu)$ as follows: let $\beta_1, \beta_2,\ldots, \beta_n$ be a basis of 
$H^1(M)$
and let $\lm_1,\ldots,\lm_n$ be  closed  one-forms with $[\lm_i] = \beta_i$
in DeRham cohomology.\footnote{When homology and cohomology  
coefficients are unspecified they are assumed to be
$\R$, so the   notation $H_1(M)$ means $H_1(M;\R)$, etc.} 
 The reader uncomfortable with 
 homology may read through this section thinking of the case
of $M=\T^n$, with angular coordinates $(x_1,\ldots, x_n)$,  and taking 
$[\lm_i]=[dx_i]$, as a basis for $ H^1(\T^n)\simeq H_1(\T^n)\simeq \R^n$. Define the
$i^{th}$  component of the rotation vector $\rho(\mu)$ as
$$
\rho_i(\mu)= \int \lm_i d\mu.
$$
Note that this integral makes sense when one looks at $\lm_i$ as inducing a function
from $TM\times\S1$ to $\R$ by first projecting $TM\times\S1$ onto $TM$, and then treating
the form as a function on $TM$ that is linear on fibers.
 The rotation vector does depend on the choice of
basis $\beta_i$, but because the  one forms are closed,  $\rho_i(\mu)$ does
not depend on the choice of representative  $\lm_i$ with $[\lm_i] = \beta_i$.
Since the rotation vector is dual to forms, 
it can be viewed as an element of $H_1(M, \R)$.
In the case $M=\T^n$, one can check that,
for $\gm$ a generic point of an {\it ergodic}  measure $\mu$, the usual rotation vector of $\gm$ coincides with that of $\mu$:
$$
\rho_i(\gm)= \lim_{b-a\to \infty}{{\tilde \gm_i}(b)-{\tilde \gm_i}(a) \over b-a}=\lim_{b-a\to \infty}{1\over b-a}\int_{d\gm|_{[a,b]}}dx_i
=\int dx_id\mu
=\rho_i(\mu)
 $$
where second equality uses the Ergodic Theorem.
 (Again, $dx_i$ is seen as a function $\PPP\to \R$). If $M$ is a general compact manifold, one can 
define the  rotation vector of a curve $\gm:\R\to M$ by
$\rho_i(\gm)=\lim_{b-a\to \infty}{1\over b-a}\int_{d\gm|_{[a,b]} }\lm_i$,
if the limit exists. As before, if $\gm(0)$
 is a generic point for an ergodic measure $\mu$,  this rotation vector does exist and coincides with that of $\mu$. 

Next we define the {\it average action}\index{average action} of a E-L invariant probability on $TM\times \S1$ by
$$
A(\mu)=\int L d\mu,
$$
which, when $\mu$ is ergodic,  we can relate to the average action 
along $\mu$-a.e. orbit $\gm$ by
$$
A(\mu)=\lim_{b-a\to \infty}{1\over b-a} \int_a^b L(\gm,\dot\gm)dt.
$$

The set of invariant probability measures, denoted $\mml$ is 
a convex set in the vector space of all measures\footnote{It is also
compact for the weak topology if, as Mather does, one compactifies $TM$.}
and   the extreme points of $\mml$ are the ergodic measures 
\cite{manerg}. Now consider the map $\mml\to H_1(M)\times \R$ given by:
$$
\mu \mapsto \lb \rho(\mu), A(\mu)\rb.
$$
This map is trivially linear and hence maps $\mml$ to a convex set
$U_L$ whose extreme points are images of extreme points of $\mml$, \ie    images of ergodic measures.
Mather shows, by taking limits of measures supported on
{long minimizers representing rational homology classes,  
that for each $\om$, there is $\mu$ such
that $\rho(\mu)=\om$ and $A(\mu)<\infty$.\footnote{The impatient reader may be tempted to
proclaim, from this fact, the existence of orbits of all rotation 
vectors. Alas, one can only deduce the rotation vector of orbits
from that of a measure when the measure is ergodic... }
Since $L$ is bounded below, the action coordinate is bounded
below on $U_L$.  Hence we can define a map
$\beta: H_1(M)\to \R$ by
$$
\beta(\om)=\inf \{ A(\mu) \mid \ \mu\in \mml, \rho(\mu)=\om\},
$$
which is bounded below and convex; the graph of $\beta$ is the boundary of $U_L$.

We say that  a probability measure $\mu\in \mml$ is a {\it minimal measure}
\index{minimal measure}
 if the point $\lb\rho(\mu), A(\mu)\rb$ is on the graph of $\beta$.
Hence, an extreme point $(\om,\beta(\om))$ of $graph(\beta)$ corresponds to 
at least one minimal ergodic measure of rotation vector $\om$.
It turns out that if $\mu$ is minimal, $\mu$-a.e. orbit lifts to a E-L minimizer in the
covering $\mb$ of $M$ whose deck transformation  group
is $H_1(M;\Z)/torsion$ (\ie the universal cover when $M=\T^n$). Conversely, if $\mu$ is an ergodic probability measure whose
support consists of $\overline M$-minimizers, then $\mu$ is a minimal measure. 

{\it Hence, each time we prove the existence
of an extreme point $(\om, \beta(\om))$, we find at least
one recurrent orbit of rotation vector $\om$ which is a $\mb$-minimizer. }

Another important property of $\beta$ is that it is superlinear, i.e ${\beta(x)\over \norm{x}}
\ra \infty$ when $\norm{x}\to \infty$. 
Since we will need the estimate later,  we motivate this in the simple case where $L=\mec$ and $\norm{\cdot}$ comes from the  Euclidean metric on the
torus.
If $\mu$ is any invariant probability measure, then
\begin{eqnarray}
\nonumber A(\mu)&=&\int L\; d\mu\geq \int \lb{\norm{\xd}^2\over 2}-\vmax \rb d\mu\\ 
&\geq &
{1\over 2}\left | \int \xd\; d\mu\right |^2 -\vmax \\  \nonumber
&= &{1\over 2}|\rho(\mu)|^2-\vmax
\end{eqnarray}
\label{superlinear}where we used the Cauchy-Schwarz inequality for the second  inequality.

The superlinearity of $\beta$ implies the existence of many 
extreme points  for  $graph(\beta)$ (although
in most cases still too few, as we will see at the end of this discussion). Indeed, $\beta$'s superlinear growth implies that
its graph cannot have flat, or linear domains going to infinity.
Any point $(\om, \beta(\om))$ is part of at least one  linear domain
of $graph(\beta)$, which we call $\lo$ (we suppress the dependence of
$c$ on $\omega$).
 Here, the index $c$ denotes the ``slope'' of the supporting hyperplane whose intersection with $U_L$ is exactly the convex and flat domain $\lo$ ($c$ can be seen as element of first cohomology). Let $\xo$ be the projection on
$H_1(M)$ of $\lo$. Then the $\xo$'s are compact and convex domains which  
``tile'' the space $H_1(M)$. Extreme points of $\xo$ are projections of extreme points of $\lo$. Hence there are infinitely many such extreme points, and infinitely many outside any compact set. Their convex hull is  $H_1(M)$, and in particular, they must
span $H_1(M)$ as a vector space. Since these extreme points are
the rotation vectors of minimal ergodic measures, we have
found that {\it there always exist at least countably many
minimal ergodic measures and at least $n=\dim H_1(M)$ of them  with distinct rotation directions. }  We will see in
Hedlund's example that this lower bound can be attained.

Finally, the generalized {\de Mather sets} \index{Mather sets (generalized)} are defined to be $\mo= support (\Mc)$, where $\Mc$ is the set of minimal measures whose rotation vector
lies in $\xo$. 
Let $\pi: TM\times \S1\to M\times \S1$ denote the
projection. Mather's main result in \cite{mather1} is the following theorem.
\begin{theorem} 
 {\bf (Mather's Lipschitz Graph Theorem)} For all $c\in H^1(M)$,
$\mo$ is a compact, non-empty subset of $TM\times \S1$. The restriction of $\pi$ to $\mo$ is injective. The inverse mapping $\pi^{-1}: \pi(\mo)\to \mo$ is Lipschitz.
\end{theorem}

In the case $M=\T^n$, Mather proves that, when they exist, KAM tori coincide with the sets $\mo$, and that they are in the closure of
these sets (see also \cite{katok} for some related results). 
In a sense, Mather's theorem generalizes Birkhoff's theorem
on invariant curves of twist maps, which says that such
curves have to be graphs (see \cite{herman}, \cite{bialypolt} for
more straightforward generalizations of this theorem}).
 The proof of the Lipschitz Graph Theorem (see \cite{mather1} or 
 \cite{mane}), which is quite involved, uses  a curve shortening argument: if curves in $\pi(\mo)$ were too close to crossing
transversally, one could ``cut corners'' and, because of recurrence,
construct a closed curve with lesser action than $A_{min}$.

\begin{remark}{\rm
An important special case is that of {\it autonomous} systems. In
this case, one can discard the time component and
view $\mo$ as a compact subset of $TM$. Then Mather's theorem implies that $\mo$ is a Lipschitz graph for the projection $\pi: TM\to M$.
To see this, suppose that  two curves $x(t)$ and $y(t)$ in $\pi(\mo)$
have $x(0)=y(s)$ for some $s$. Mather's theorem rules out immediately the possibility that $s$ is an integer, unless $x=y$ is a periodic orbit. For a general $s$, consider the curve $z(t)=y(t-s)$.
Then, $\dot z(t)= \dot y (t-s)$ and, by time-invariance of the Lagrangian, $(z(t),\dot z(t))$ is a solution of the E-L flow. It
has same average action and rotation vector as $(y,\dot y)$ and hence it is also in $\mo$. But then $z(0)=x(0)$ is impossible, by
Mather's theorem,
unless $\dot z(0)=\dot y(s)= \xd(0)$.

In the realm of twist maps, one can also deduce from Mather's theory the existence
of many invariant sets that are graphs over the base and are made of minimizers (see also \cite{katok}).}
\end{remark}

\subsection{Examples}

One would hope that $\beta$ is, in general, strictly
convex, \ie  each point on $graph(\beta)$ is an extreme point.
This is true when $M=\S1$, and Mather shows in \cite{mather2}
 how his Lipschitz Graph
Theorem implies the classical Aubry-Mather Theorem, by taking
a E-L flow that suspends the twist map. The fact
that $\mo$ is a graph nicely translates into the fact that orbits in an
Aubry -Mather set are well ordered. 

The graph of $\beta$ is also strictly convex  when $L$ is a Riemannian metric on $\T^2$. This was known by Hedlund \cite{hedlund}  in the 30's, albeit in a different language. It  will also be an easy consequence of Theorem \ref{t2}. In Section 6, we will show that this may not be true if one adds a  potential term  to the metric in the Lagrangian.

Here we briefly describe three other counter-examples to the 
strict convexity of $\beta$.

\begin{example} {\bf (Ma\~ne \cite{mane})} 
\label{maneg} {\rm Take $L: T\T^2\to \R$,
given by $L(x,\xd)= \norm{\xd-X}^2$ where $X$ is a vector field on $\T^2$.
The   integral curves $x$ of $X$ are
 automatically E-L minimizers
since $L\equiv 0$ on these curves. Ma\~ne chooses the vector field $X$ to be a (constant) vector field of irrational slope multiplied  by a carefully chosen
function on the torus which is zero at exactly  one point $q$. 
 The integral flow of $X$ has the  rest point $q(t)=q$, and all the other solutions
 are dense on the torus. The flow of $X$ (and its lift to $T\T^2$ by 
the differential) has 
 exactly two ergodic measures:
one is the  Dirac measure 
supported on $(q,0)$, with zero rotation vector, the other is equivalent to the
 Lebesgue measure on $\T^2$ and has nonzero rotation vector, say $\om$ (see \cite{handel}, for more details
on this $\T^2$ flow). 
 Ma\~ne checks that $\beta^{-1}(0)$ (trivially always an $\xo$)
 is the interval $\{\lm\; \om \mid\ \lm \in [0,1]\} $,  and that  no
 ergodic measure has a rotation vectors strictly inside this interval. Thus the
Mather set $M_0$ is the union of the supports of the two measures.}
\end{example}

In the context considered here,
this example is a little unsatisfactory  because it is not
a mechanical Lagrangian. We will give an example in Section 6 of a mechanical 
Lagrangian on $\T^2$ which displays a similar phenomena. The next
two examples involve metrics, on the three-torus and on the surface
of genus two, respectively.

\begin{example} {\bf (Hedlund-Bangert)} 
\label{hedlundeg} {\rm Consider in $\R^3$ the three nonintersecting lines given by the $x$-axis, the $y$-axis translated by $(0,0,1/2)$ and the $z$-axis translated by $(1/2,1/2,0)$. Construct a $\Z^3$-lattice
of nonintersecting axes by translating each one of these by all integer vectors.
 Take a metric in $\R^3$ which
is the Euclidean metric everywhere except in small, nonintersecting
tubes around each of the axes in the lattice. In these tubes, multiply
the Euclidean metric by a function $\lm$ which is  1 on the boundary and attains its (arbitrarily small)
minimum along the points in the center of the tubes, \ie   at the axes
of the lattice.  Because the construction is $\Z^3$ periodic, this metric induces  a Riemannian metric
on $\T^3$. One can show (\cite{bangert}), if $\lm$ is taken sufficiently small, that a minimal geodesic (which is a E-L minimizer in our context) can make at most three jumps between tubes. In particular, a recurrent E-L minimizer has to be one of the three disjoint periodic orbits which are the projection of
 the axes of the lattice. Thus there are only three rotation directions that
minimizers can take in this example, or six if one counts positive and
negative orientations. In terms of Mather's theory,
the level sets of the function $\beta$ are octahedrons with vertices
$(\pm a, 0,0), (0,\pm a, 0), (0,0,\pm a)$ (we assume here that $\lm$
is the same around each of the tubes). Since we are in the case of a metric,  one can check that $\beta$ is quadratic when restricted to a line
through the origin (a minimizer of rotation vector $a\om$ is a reparametrization
of a minimizer of rotation $\om$). Hence a set $\lo$ is either a face,  an edge
or a vertex of some level set $\{\beta=b\}$, and the corresponding $\mo$ is,
 respectively, the union of three, two (parameterized at same speed) or one of the minimal periodic orbits one gets by projecting the disjoint axes. Note
that, instead of the function $\beta$ of Mather, Bangert uses the {\it stable
norm}. Mather's function $\beta$ is a generalization of that norm.}
\end{example}

It is important to note  that the nonexistence of minimizers of a certain 
rotation vector $\om$ does not mean that there are no  orbits of the E-L flow that  have  rotation vector $\om$.
For example, Mark Levi has shown the existence of
 orbits of all rotation vectors
in the Hedlund example (personal communication). In addition,
 in our torus example (Section 6), it is
easy to see that there are (nonminimizing) orbits with rotation vectors in
the excluded interval.
The next example was brought to our attention by
A. Fathi.

\begin{example}
\label{fig8} {\rm Take the metric of  constant negative curvature on 
the surface of genus 2 (the two-holed torus) which has a thin neck between the
two holes (see Figure 1).
 In this case, the notion of minimizing is just that of least
length using the hyperbolic metric.  With $a$ and $b$ as shown, the minimal measure
for the homology class $a + b$ will be a linear combination of the ergodic
measures supported on $\Gamma_a$ and $\Gamma_b$, where
$\Gamma_a$ and $\Gamma_b$ are the closed geodesics in the homotopy classes
of $a$ and $b$, respectively. This is because any closed curve that crosses the
neck will be longer than the sum of the lengths of $\Gamma_a$ and $\Gamma_b$. Hence $(a+b, \beta(a+b))$ cannot be an extreme point of
$graph (\beta)$. }
\end{example}

\begin{figure}
\centerline{\psfig{figure=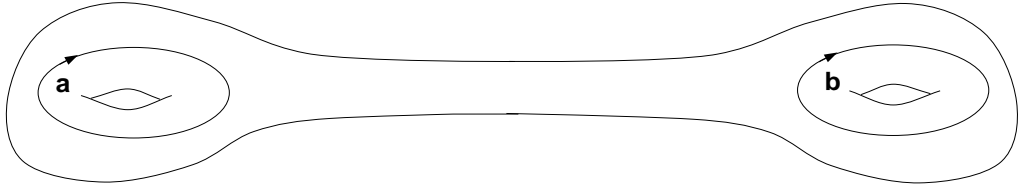,width=.9\hsize}}
\caption{The genus two surface in Example 5}
\end{figure}

This example illustrates the remarks of the introduction, namely,
on hyperbolic manifolds the notion of $\mb$ minimizers 
is not the correct one if
one wants to show that all the dynamics of the geodesic flow
of a metric of constant negative curvature are preserved  under
global perturbation. The right notion of minimality comes from the observation
that there is a curve of least length in the {\it homotopy} class of $ab$. 
This notion is generalized to an asymptotic invariant
on hyperbolic manifolds that associates an ergodic measure of a given
dynamical system with its ``rotation measure''. The rotation measure is
the ergodic measure of the hyperbolic geodesic flow whose dynamics
best mirrors that of the given ergodic measure (see \cite{philrm}).

\section{ Autonomous Lagrangian  Systems on the 2-Torus}

We now give a concrete application of Mather's theory to
autonomous Lagrangians on the two torus. If the Lagrangian is
mechanical, then for each $E > \vmax$ we get a Jacobi metric.
Hedlund showed (\cite{hedlund}) that for each Riemannian metric on the
two torus there are minimizing geodesics in all directions. Applying
this to the Jacobi metric we see that there will be minimizers in all
directions on each energy level. However, the parameterization of these
minimizers as solutions of the E-L flow is different than as Jacobi
geodesics. Thus in applying Hedlund's result we must examine how the
minimizers on the various levels fit together in terms of their rotation
vectors. This can be done nicely in the framework of
Mather's theory which also allows us to consider a
wider class of Lagrangians.

Recall
that the sets $\xo$ are projections on the rotation vector space of flat domains in $graph (\beta)$.

\begin{theorem} 
\label{t2}
For any autonomous Lagrangian system on $\T ^2$ which satisfies Mather's hypothesis,  
a set $\xo$ is either  a finite interval in  a line through the
origin
 or a point. If $\xo$ is an interval, either it contains 0,
 or it is supported by a line of rational slope. In addition:
\begin{description}
\item{(i)}  There are minimal  ergodic measures 
(and hence recurrent E-L minimizers) of all rotation directions.
\item{(ii)}  If $v$ has irrational direction, then there is a positive number
 $K(v)$ such that  $\norm{v}> K(v)$ implies that
 there is a minimal ergodic measure (and hence a recurrent E-L minimizer)
with rotation vector $v$.
\item{(iii)} In each rational direction, represented by a vector $v$, 
there is a sequence of positive numbers $\lm_n\to \infty$ and 
a sequence  $\gm_n$ of periodic E-L minimizers such that $\gm_n$
has rotation vector $\lm_nv$.
\item{(iv)} If the Lagrangian is of the form $L(x,\xd)= {1\over 2}\norm{\xd}^2 -V(x)$,  where $\norm{\cdot}$ is the Euclidean metric,\footnote{Clearly, similar results can be obtained for other metrics} then the origin is not 
in the interior of any $\xo$ and the $\xo$ which have the origin as endpoint  have  length less than 
$2\sqrt{V_{max}-V_{min}}$. In particular, $K(v)$ in (ii) can be taken to be
$2\sqrt{\vmax- V_{min}}$.

\item{(v)} The support of  minimal measures can
be either a point, a closed curve, the suspension of a Denjoy example (a lamination), or the whole torus.
\end{description}
\end{theorem}

\Bproof By Mather's Lipschitz Graph Theorem and the remark after it, any set
$\mo$
projects injectively to $\T^2$. It thus yields a flow on compact set $X\subset \T^2$.
  Most of the theorem is a consequence of the fact that there are
severe restrictions on what $X$ and its flow can be.
We review some of these well known results. They follow
easily from, for example, the techniques of \cite{franksmis} (
\cf\ \cite{walsh}).

If $X\subset \T^2$ is a compact set and there is a continuous flow
$\phi_t$ defined on $X$, then the support of a $\phi_t$-invariant
ergodic measure on $X$ can only be
(1) a point, (2) a closed orbit,
(3) the whole torus, or (4)  homeomorphic
to the suspension of the minimal set in a homeomorphism of the circle
 that is a
Denjoy counter example (and it must be embedded in $\T^2$ in the obvious way).

In cases (1) and (2) there is a unique invariant probability
measure on the set and so the rotation vector exists and
is equal for all points.  Case (1) happens only if
the rotation vector is $0$. Case (2) can happen only if the rotation
vector is zero for homotopically trivial curves and
a nonzero $v\in\Q^2$,
if the curve is homotopically nontrivial.
Cases (3) and (4) are more subtle. The set $X$
can contain a fixed point, but in
any event the rotation vectors of ergodic measures supported on
$X$ can take on at most two values, zero  and a  vector with irrational slope.

Of particular importance in what follows are the following three
consequences. First, there cannot be two invariant measures in
$X$ whose rotation vectors are in different directions,  \ie  
${v\over \norm{v}} \not ={w \over \norm{w}}$. Second, if an
ergodic measure has a nonzero rotation vector
$v\in \Q^2$, then it is supported on a closed orbit. Finally, $X$
cannot support two ergodic measures with nonzero rotation vectors
that have the same irrational slope.

To prove the first paragraph of the theorem, recall that by Mather's theory,
 the endpoints of each $\xo$ corresponds to a minimal
 ergodic measure with that rotation vector. 
Next observe that  any $\xo$ is a finite interval on a line through the 
origin or a point, because as  noted in the previous paragraph, no $\mo$ can 
 contain orbits with different nonzero
 rotation directions and each $\mo$ can contain ergodic measures that yield
at most two rotation vectors. (The fact that   $\xo$ is finite also follows
from the superlinearity of $\beta$.) If we 
 fix an $\omega$ with irrational slope, then if both  endpoints of the  interval $\xo$
are nonzero, they correspond to at least two distinct ergodic measures with  two 
irrational rotation
vectors of same direction but different lengths, another impossibility.

 To prove (i), take any $\omega\in\R^2$. Using what we have just proved, any set
$\xo$ to which $\om$ belongs is a finite interval with endpoints in the same direction as $\omega$. 
These endpoints are the rotation vectors of two minimal ergodic measures. 
For (ii),  note that if $\omega$ has irrational direction,
 and if $\norm{\omega}$ is large enough, we have just seen
 that $\xo$ can only be a point. Hence there is a minimal
ergodic measure of rotation vector $\omega$.
 As for (iii), note that the previous paragraph implies, 
in each rational rotation direction $v$,
 the existence of sequences $\lm_n\to \infty$ and minimal
 ergodic measures $\mu_n$ with $\rho(\mu_n)=\lm_nv$.
 As noted at the beginning of the
proof, such an ergodic measure  necessarily comes from a periodic orbit.

We  now prove part (iv). 
We first show that $\beta(0)$ is an isolated global minimum,
thus proving that $0$ cannot be  in the interior of any $\xo$.
In the case of autonomous mechanical systems, it is easy
to see that a measure that gives the absolute minimum for
$A$ is the Dirac measure on  a fixed point  $x$ at which $V(x)=\vmax$.
 This measure has rotation vector $0$ and action $\vmax$. 

Recall, moreover, that from Formula (2) in Section 4,  
 we have the estimate
 $\beta(\omega)\geq {1\over 2}\norm{\om}^2 -\vmax$ for such systems.
This implies
that the graph of $\beta$ above an $\xo$  that contains zero 
 cannot be horizontal, and that $(0,\vmax)$ must be an isolated minimum of  $\beta$. 

Finally, we prove our estimate for the length of $\xo$'s
which have  0 as one of their endpoints.
On an energy level $E>\vmax$ the Jacobi metric 
$\sqrt{E-V)}ds$ is a 
Riemannian metric on $\T^2$. It is not hard to see that
the Jacobi minimizers and the E-L minimizers for $L_E={1\over 2}
(E-V)\norm{\xd}^2$  form the same set of curves on the torus
 (See \cite{milnor}, page 70). 
$L_E$ obviously satisfies Mather's hypothesis and so, by Part (i), there are Jacobi minimizers of all rotation directions 
(we could also have invoked Hedlund \cite{hedlund}, who proved this in the 30's).
 By Proposition \ref{jel}, 
 these Jacobi minimizers are 
also E-L minimizers for $L=\mec$. In terms of minimal ergodic
measures, after a time change  we have found minimal ergodic measures of all
rotation  direction  in each energy level $E>\vmax$.

Let  $\mu$ be such a measure. Using  Lemma \ref{disac},  we have:
$$
{1\over 2}\left(\rho(\mu)\right)^2 - \vmax \leq \int L d\mu
\leq \int E-2V_{min} d\mu= E-2V_{min}
$$
so $\rho(\mu)\leq \sqrt{2(E-2V_{min}+\vmax)}$. This is true  for the measures obtained  from
 energy levels $E$ arbitrarily close to $\vmax$. 
 Thus we have shown that $\norm{\om}\leq 2\sqrt{\vmax-V_{min}}$. 

Statement $(v)$ is an immediate consequence of the facts 
at the beginning of the proof.
 \Eproof

\section{ Speed Defect: an Example in the  2-Torus}

\subsection{The Main Example}

In this section, we exhibit an autonomous Lagrangian system on the
torus $\T ^2$ which fails to have recurrent minimizers of all rotation vectors in a prescribed rational direction.  This is a different kind of counterexample to a generalization of the Aubry--Mather theorem  than  Hedlund's example (one can also derive counterexamples in the realm of symplectic twist maps using examples in \cite{arnaud}).
In terms of Mather's theory, we  find in our system  an $\xo$ in a
rational direction which
contains an interval, and hence the function $\beta$ of Mather
is not strictly convex in this case, confirming the worst predictions
of Theorem \ref{t2}.   

 Note that the Lagrangian in such an example could not
be a Riemannian metric, since it is easy to show that Mather's function 
$\beta$ 
for such systems is quadratic along lines through the origin of $\R^2$. Hence
Theorem 5 above implies that, in the case of metrics,  $\beta$  is strictly convex.

We first outline the features of our example. We will find, in a given
energy level $(h=0)$, exactly two closed minimizers $\gm_1, \gm_2$ in the 
homology class $(0,1)\in \Z^2=H_1(\T ^2;\Z)$
but with different rotation vector $ (0,\rho_1), (0,\rho_2)$, \ie
 different average speeds. We  then
show that any periodic minimizer in that same homology class but in a different
energy level goes strictly faster for higher energy, and strictly slower
for lower energy than these two minimizers, leaving the segment between
$(\rho_1,0)$ and $(\rho_2, 0)$ empty of recurrent minimizers.

On the torus $  \T^2= \R^2/{2\pi \Z^2}$, consider a Lagrangian of the form:

\begin{eqnarray}
L(x,\xd)&=& {\alx\over 2 }\norm{\xd}^2 - V(x_1)\\
V(x_1) &=&-\frac{2+\cos 2x_1}{2+\sin x_1}\\
\alpha (x_1) &=& {2+\sin x_1},
\end{eqnarray}
where $\norm{\xd}$ denotes the usual, Euclidean norm on $\T ^2$.
The corresponding Hamiltonian  is $H(x, p)= {\norm{p}^2\over 2\alx } + V(x_1)$.

If there is a minimal ergodic measure in the homology direction $(0,1)$,
 then
as noted at the beginning of the proof of Theorem \ref{t2}, it must be
supported on a periodic orbit.
Furthermore, we claim that if $\gm$ is
 a closed curve which is the projection of the
 support of an ergodic, minimal measure in
the $(0,1)$ direction, it has to be of the form $x_1= constant$. Indeed, since
$H(x,p)$  is independent of $x_2$,  the Hamiltonian flow is invariant
under translations in the $x_2$ direction. Hence $\gm$ is actually part of
a one parameter family of periodic orbits that are translates of one another
in the $x_2$ direction. Obviously, all the orbits in this family have same
rotation vector and action. Since $\gm$ is the support of
a minimal measure, all its $x_2$ translates are as well, and they  all belong to the same 
$\xo$. Because $\gm$ is in the class $(0,1)$, it must have (at least 2) points
at which $p_1=0$. If it is not so at all points, \ie   if $x_1$ is not
constant along $\gm$, it is easy to see that $\gm$ will intersect transversally
its $x_2$ translates, contradicting (by Mather's theorem) the fact that these curves belong to the 
projection of the same $\xo$.

Our task is now to find all the possible periodic orbits of the Hamiltonian
flow that project to curves $x_1=constant$ and to compute their rotation
vector and action.
Such a curve occurs when $0=\dot x_1={p_1\over\alx}$ and hence (by setting
$p_1(0)=0$), whenever
$$
0\equiv \dot p_1(t)= -{\del H\over\del x_1}={\alpha'\over 2\alpha^2}\norm{p}^2 -V'
$$
Simplifying this equation by fixing an energy level $H=E$ and using $\norm{p}^2= 2\alpha(E-V)$, one
gets:

$$
0=\cos x_1 (2+\sin x_1)(4\sin x_1 -E)
$$
which has the solutions:

\begin{eqnarray} 
x_1&=&{\pi\over 2},  \\
 x_1&=&{3\pi\over 2},  \\
4\sin x_1&=& E, 
\end{eqnarray}
with some restrictions on what the range of $E$ is in each case, which we will deal with later.

We now find the curves that the corresponding ergodic measures trace in the 
$(\norm{\omega}=\rho, A)$- plane, and show that the third curve is always
in the union of the epigraphs of the two first ones, and hence cannot
be a minimal measure (the epigraph of a function is the set
of points above its graph) .

In the energy level $\{H=E\}$,  the Lagrangian is given by
$L=E-2V$. Since it is constant along the curves $x_1=constant$,
the average action of the corresponding  ergodic (probability) measure $\mu$ is 
$$
A(\mu)=\int Ld\mu=L\int d\mu=L
$$
From $L=E-2V$, we can also get the rotation vector
 $(0,\rho)$ of an $x_1=constant$ curve:
$$
L+2V={\alx\over 2} \norm{\xd}^2+V(x_1)=E\implies
 \norm{\xd}^2= \frac{2(E-V(x_1))}{\alx}
$$
In particular, the speed at which a curve $x_1= constant$ is traversed
by the E-L 
flow is constant. The rotation vector of such a curve is:

$$
(0,\rho(E,x_1))= \lb 0, {\norm{\xd}\over 2\pi}\rb=
 \lb 0, {1\over 2\pi}\sqrt{2(E-V)\over \alpha}\rb
 $$

We found it more convenient to first look at the curves $(c,A)$ formed
by these measures, where 
$$
c \stackrel{\rm def}{=} {(E-V)\over \alpha} = 2\pi^2\rho^2.
$$
We denote by $(c_i(E),A_i(E))$ the curve corresponding to
 the $i^{th}$ family of solutions, $i=1,2,3$.

$\bullet$ If $x_1={\pi\over 2}$, $A_1(E)=L= E-2V({\pi\over 2})=E+{2\over 3}$,
whereas $c_1(E)={ E-V({\pi\over 2})\over \alpha({\pi\over 2})}={1\over 3}(E+{1\over 3})$. Hence:
$$
A_1(c)=3c+{1\over 3}.
$$
One can verify that the curve $x_1={\pi\over 2}$ is indeed a 
solution for $E$ in $[-1/3,\infty)$, which gives $c \in [0,\infty)$.

$\bullet$ If $x_1={3\pi \over 2}$, $A_2(E)= E+2$, $c_2(E)= E+1$ and hence
$$
A_2(c)= c+1.
$$ 
One can verify that $x_1={3\pi\over 2}$ is a solution for $E$ in $[-1,\infty)$,
which corresponds to $c\in [0,\infty)$.

$\bullet$ If $E=4\sin x_1$, we get, by replacing $\sin x_1=E/4$ and
$\cos 2x_1= 1-2\sin^2x_1= 1-2(E/4)^2$:

\begin{eqnarray*}
c_3(E)&=&2 {E^2+16E+24\over (E+8)^2}\\
A_3(E)&=& 8{(E+3)\over (E+8)}.
\end{eqnarray*}

We will be content to say that the curves $E=4\sin x_1$ correspond to solutions
for $E$ in some interval included in $[-2,4]$ (one can check that $-2<V_{min}<-1$), which we do not need to find, as we will see.

 
Instead of trying to find $A_3(c)$ explicitly, we compute $A_3(E)-(c_3(E)+1)=
{5E^2-8E+80\over (E+8)^2}$ which is positive for all $E$, thus showing that the 
curve $(c, A_3(c))$ is above $(c, A_1(c))$. The
relation is trivially the same for the corresponding $(\rho, A(\rho))$
curves and hence the measures corresponding to the curves $E=4\sin x_1$ can
never be minimal.

\begin{figure}
\centerline{\psfig{figure=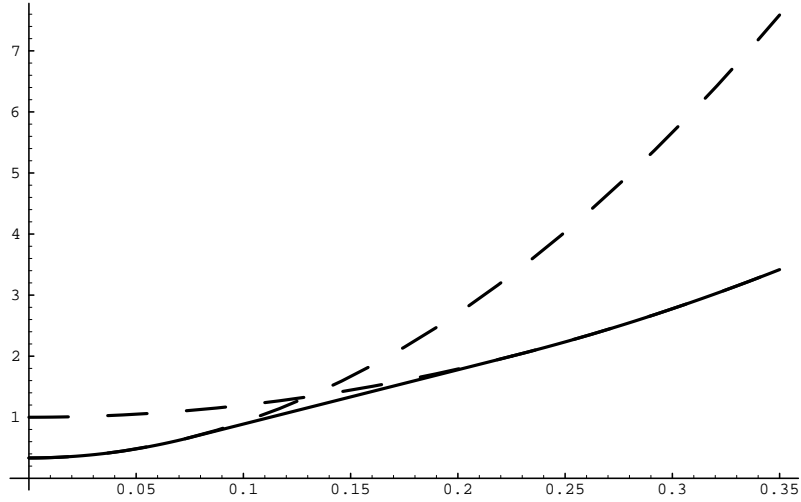,width=.7\hsize}}
\caption{The graph of $\beta$ in the $(0,1)$ direction}
\end{figure}

 Replacing $c$ by $2\pi^2 \rho^2$, we get:
\begin{eqnarray*}
A_1(\rho)&=&6\pi^2\rho^2+{1\over 3} \\
A_2(\rho)&=& 2\pi^2 \rho^2 +1,
\end{eqnarray*}
 the graphs of which are two parabolas crossing at $\rho= {2\pi^2\over 3}$,
with $A_1(\rho)\leq A_2(\rho) $ for $\rho \in [0, {2\pi^2\over 3}]$, $A_2<A_1$ after that. Note that the union of the epigraphs of these curves
is not convex. Remember that the graph of  Mather's function $\beta$ is the boundary of the convex
hull of the points $(\rho (\mu),A(\mu))$ for all possible ergodic
measures $\mu$ .
To form the lower boundary of the convex hull for the epigraphs of $(\rho, A_i(\rho))$, $i=1,2$ (which, as we showed, necessarily contains the epigraph
of $(\rho, A_3(\rho))$, we look for the line tangent to both these graphs.
It is easy to check that it joins the point $Q_1=({1\over 3\sqrt 2 \pi}, 2/3)$
of the graph of $A_1$ to the point $Q_2=({1\over \sqrt 2 \pi}, 2)$
of the graph of $A_2$. Hence, any measure corresponding to points on
the graphs of $A_1$ or $A_2$ above the segment of line between $Q_1$ and $Q_2$
fails to be  minimal. Instead, if $\mu_i$ are the measures corresponding to 
$Q_i$, $i=1,2$, the minimal measure with rotation vector in the interval
$I=[ {1\over 3\sqrt 2 \pi},{1\over \sqrt 2 \pi}]$ will be a convex combination of
$\mu_1$ and $\mu_2$. The set $\xo$ for any $\omega=(0,\rho)$ with $\rho\in I$
is thus an interval, and the corresponding Mather set $\mo$ is the union of
the support of $\mu_1$ and $\mu_2$, \ie   2 loops, with
rotation vectors $(0, {1\over 3\sqrt 2 \pi})$ and $(0, {1\over \sqrt 2 \pi})$. Since any recurrent E-L minimizer with rotation number $\omega$ belongs to the Mather
set $\mo$, it
follows that there can't be any recurrent minimizer of rotation vector
$(0,\rho)$ with $\rho$ in the interior of $I$. This is the gap we announced.

\begin{remark} 
{\rm It is easy to check that both $\mu_1, \mu_2$ have their support in $\{H=0\}$,
and that $(\rho, A_1(\rho))$ for $ \rho \in [0,{1\over 3\sqrt 2 \pi}]$ corresponds
to measures supported in $E\leq 0$ whereas $(\rho, A_2(\rho)), \rho\geq {1\over \sqrt 2 \pi}$ corresponds to $E\geq 0$. 

It is instructive to think of the example in terms of the Jacobi metric.
For $E=0$, this metric is $\sqrt{2 \cos(2 x_1)} \norm{\dot x}$, which
has exactly two minimizers (with respect to the metric) in the homology
class of $(1,0)$, namely, $x_1 = \pi/2$ and $x_1 = 3 \pi/2$. Of necessity these
closed curves have the same Jacobi length, {\it but} when they are parameterized
as solutions to the E-L equations, they have different speeds and thus
different rotation vectors.

As the energy $E$ is swept through $0$ there is a ``minimality exchange''.
Namely, when $E\leq 0$, $x_1 = \pi/2$ is the Jacobi minimizer and corresponds to
a minimal measure for the E-L flow. When $E\geq 0$, $x_1=3\pi/2$
has these properties. At $E=0$ they share these properties, but as noted above,
they have different rotation vectors for the E-L flow.}
\end{remark}

\subsection{A E-L minimizer which is not a Jacobi minimizer}

Using the example above we construct a E-L minimizer that is not
a Jacobi minimizer on its energy level, thus showing that the converse
of Proposition \ref{jel} in Section 2.4 is false.

Let $\tae$ be the parameterization of the curve $x_1\equiv\pi/2$ with
the property that $d\tae$ is a trajectory of the E-L flow with $H = E \in
[1/3, \infty)$. Similarly, $\tbe$ corresponds to the curve $x_1\equiv
3\pi/2$ for $E\in[-1,\infty)$. We shall use $\alpha$ and $\beta$ to refer
to these curves without a specific parameterization.

 If $\mu_E$ and $\eta_E$ are the ergodic
invariant probability measures supported on $d\ae$ and $d\be$,
respectively, then we have that for $E<0$ and $E>0$, respectively, 
$\mu_E$ and $\eta_E$ are minimal measures. Further, $\mu_0$ and
$\eta_0$ are the minimal measures corresponding to the rotation vectors
$\rho_1 = (0, \frac{1}{3 \sqrt{2}\pi})$ and
 $\rho_2 = (0, \frac{1}{\sqrt{2}\pi})$ and
for any $\rho = (0, p)$ with $\frac{1}{3\sqrt{2}\pi} < p <
\frac{1}{\sqrt{2}\pi}$,
the minimal measure is a linear combination of $\mu_0$ and $\eta_0$.
In particular, if $p_3 = \frac{2}{3\sqrt{2}\pi}$ and $\rho_3 = (0,p_3) = {\rho_1
+ \rho_2\over 2}$, then the minimal measure with rotation vector
$\rho_3$ is $\frac{\mu_0 + \eta_0}{2}$.

Now for each $N\in\N$, let $\gmn$ be the E-L minimizer satisfying
$\gmn(-N/2) = (\pi, -N p_3/2)$ and 
$\gmn(N/2) = (\pi, N p_3/2)$, and so $\rho(\gmn) = \rho_3$. If we let
$\sigma_N$ be the probability measure uniformly distributed with respect
to time on $d\gmn$, then the results in section 2 of \cite{mather1} says
that as $N\ra\infty$, $\sigma_N$ will converge weakly to the minimal
measure with rotation vector $\rho_3$ and so $\sigma_N\ra\frac{\mu_0 +
\eta_0}{2}$. In particular, by choosing $N$ large, we can make $H(\gmn,
\dot\gmn)$ arbitrarily close to zero, and we can insure that $\gmn$ spends
about half its time near $\alpha_0$ and the other half near $\beta_0$.
 For each $N$, let $\ell_N$ be the
path that is the union of 
$$[(\pi/2, \frac{Np_3}{2}), (\pi, \frac{Np_3}{2})],
 [(\pi/2, \frac{Np_3}{2}), (\pi/2, \frac{-Np_3}{2})], \mbox{\rm and\ }
 [(\pi/2, \frac{-Np_3}{2}), (\pi, \frac{-Np_3}{2})],
$$
where we
use square brackets to indicate the straight segment connecting two points
in the plane. Thus $\ell_N$ goes directly left from the top endpoint of
$\gmn$, then travels down $\alpha$ and then directly across to the
bottom endpoint of $\gmn$. 

Now we assume that for all $N$, $\gmn$ is a Jacobi minimizer for its
energy $E_N$ and obtain a contradiction. For $N$ large, we have that 
$|E_N|$ is small. Now if $E_N \leq 0$, since $\ae$ is the minimizer
for those $E$, it is clear that $\ell_N$ has lesser Jacobi length than
$\gmn$ for energy $E_N$ because $\gmn$ must spend half of its time near
$\beta$. By considering a curve analogous to $\ell_N$ but going to the
right and then down $\beta$, we find that $E_N\geq 0$ is impossible also.

\begin{remark}
{\rm As noted in Section 3.3
and Example 4 in Section 4.2, when one attempts to take the limit
of minimizers $\gmn$ without adequate control, it can happen that
the $\gmn$ converge to something that is not in the same direction as
each $\gmn$. The example of this section
illustrates another mechanism by which the
limit argument can fail. It can happen that the $\gmn$ converge to
something in the correct direction and  the speed of the limiting
measures average to the correct speed, but the various ergodic
components of the limit measure have either
greater or lesser speed than the desired one. This gives rise to
the speed gap one has to allow in Theorem 1 and Theorem 5 .}
\end{remark}

\section*{Acknowledgments} The authors would like to thank Robert MacKay for
providing inspiration for  this work.
 We would also like to thank the IMS at Stony Brook and all of its members for 
their support while a good part of this work was done.

\section*{Note Added in Proof} Just prior to the publication of this
paper we became aware of the elegant paper ``On minimizing measures of
the action of autonomous Lagrangians'', by M. J. Dias Carneiro 
({\it Nonlinearity}, {\bf 8}, 1077--1085, 1995)  which 
contains much that is relevant to Sections 5 and 6 of this paper.
In particular, Dias Carneiro proves that in the autonomous
case,  the function $\beta$ (defined in Section  4.1)
has a derivative in radial directions, 
calculates the derivative,  and as a corollary, 
 shows that  Mather sets are always supported on a single energy level.
These results simplify the required calculations in the  example 
of Section 6.
In addition, Dias Carneiro makes
a remark that contains part of Theorem 5 and gives  a result
connecting the minimal measures of the Jacobi metric with those of
the E-L flow.

In addition, R. Iturriaga showed us notes of \Mane\  which contain
11 beautiful  theorems about autonomous Lagrangian systems. His  students are 
providing the proofs of these results. 

We dedicate this paper to the memory of Ricardo \Mane.

\end{document}